\documentclass[a4paper,10pt]{amsart}
\usepackage{amssymb,amsmath,amsfonts,amsthm}
\usepackage[dvips]{graphicx}
\usepackage{psfrag}

\theoremstyle{plain}
\newtheorem{main}{Theorem}

\newtheorem{theorem}{Theorem}[section]
\newtheorem{lemma}[theorem]{Lemma}
\newtheorem{proposition}[theorem]{Proposition}
\newtheorem{corollary}[theorem]{Corollary}
\theoremstyle{remark}

\newcommand{\quand}{\quad\text{and}\quad}
\newcommand{\Leb}{\operatorname{vol}}
\newcommand{\per}{\operatorname{per}}

           \def\ea{\end{array}}
          \def\ec{\end{center}}
     \def\ed{\end{description}}
        \def\ee{\end{equation}}
       \def\eea{\end{eqnarray}}
     \def\eeaa{\end{eqnarray*}}
 \def\et{\end{thebibliography}}

\def\Orb{{\rm Orb}}
\def\Diff{{\rm Diff}}

\def\supp{\operatorname{supp}}

\def\cC{{\mathcal C}}
\def\cO{{\mathcal O}}

\def\cU{{\mathcal U}}
\def\cV{{\mathcal V}}

\def\cF{{\mathcal F}}
\def\cM{{\mathcal M}}

\def\cS{{\mathcal S}}
\def\cW{{\mathcal W}}

\def\id{\operatorname{id}}

\def\vep{\varepsilon}

\def\TT{{\mathbb T}}
\def\RR{{\mathbb R}}

\def\NN{{\mathbb N}}

\title[Maps with mostly contracting center]{Geometric and measure-theoretical structures of maps with mostly contracting center}
\author{Dmitry Dolgopyat, Marcelo Viana and Jiagang Yang}
\date{\today}
\thanks{M.V. and J.Y. were partially supported by CNPq, FAPERJ, and PRONEX}

\address{4417 Mathematics Building, University of Maryland, College Park, MD 20742, USA}
\email{dmitry@math.umd.edu}

\address{IMPA, Est. D. Castorina 110, 22460-320 Rio de Janeiro, Brazil}
\email{viana\@@impa.br}

\address{Departamento de Geometria, Instituto de Matem\'atica e Estat\'\i stica, Universidade Federal Fluminense, Niter\'oi, Brazil}
\email{yangjg\@@impa.br}

\begin{document}

\begin{abstract}
We show that every diffeomorphism with mostly contracting center direction exhibits a geometric-combinatorial
structure, which we call \emph{skeleton}, that determines the number, basins and supports of the physical measures.
Furthermore, the skeleton allows us to describe how the physical measure bifurcate as the diffeomorphism changes.
In particular, we use this to construct examples with any given number of physical measures, with basins densely
intermingled, and to analyse how these measures collapse into each other - through explosions of their basins -
as the dynamics varies. This theory also allows us to prove that, in the absence of collapses, the basins are
continuous functions of the diffeomorphism.
\end{abstract}

\maketitle

\setcounter{tocdepth}{1} \tableofcontents


\section{Introduction}

The notion of mostly contracting center refers to partially hyperbolic diffeomorphisms and
means, roughly, that all Lyapunov exponents along the invariant center bundle are negative.
It was introduced by Bonatti, Viana~\cite{BoV00}
as a more or less technical condition that ensured existence and finiteness of physical measures.
Since then, it became clear that maps with mostly contracting center have several distinctive
features, that justify their study as a separate class of systems.

For instance, Andersson~\cite{An10} proved that they form an open set in the space
of $C^{1+\vep}$ diffeomorphisms, and that the physical measures vary continuously on an open
and dense subset. Castro~\cite{Cas02,Cas04} and Dolgopyat~\cite{Dol00} studied the mixing
properties of such systems. Moreover, Dolgopyat~\cite{Dol04a} obtained several limit theorems
in a similar context. In addition, Melbourne, Matthew~\cite{MN05} proved an almost sure invariance
principle (a strong version of the central limit theorem) for a class of maps that includes
some partially hyperbolic diffeomorphisms with mostly contracting center.
Burns, Dolgopyat, Pesin~\cite{BDP02} studied maps with mostly contracting center in the
volume preserving setting, obtaining several interesting results about ergodic components,
stable ergodicity, and other aspects of the dynamics.
Moreover, Burns, Dolgopyat, Pesin, Pollicott~\cite{BDPP08} studied stable ergodicity of
Gibbs $u$-states, in the general (non-volume preserving) setting.

Before all that, Kan~\cite{Ka94} exhibited a whole open set of maps on the cylinder
with two physical measures whose basins are both dense in the ambient space.
His construction was extended by Ilyashenko, Kleptsyn, Saltykov~\cite{IKS08}.
See also~\cite[\S~11.1.1]{Beyond}.
As it turns out, these maps have mostly contracting center. This construction can also be
carried out in manifolds without boundary, but then it is not clear whether coexistence
of physical measures can still be a robust phenomenon. This is among the questions we
aim to answer in this paper: we find negative answers in some situations.

Systems with mostly contracting center have been found by several other authors.
Let us mention, among others:
Ma\~n\'e's~\cite{Man78} examples of robustly transitive diffeomorphisms that are not hyperbolic
(see also~\cite{BoV00} and \cite[\S~7.1.2]{Beyond});
Dolgopyat's~\cite{Dol04b} volume preserving perturbations of time one maps of Anosov flows;
volume preserving diffeomorphisms with negative center Lyapunov exponents and minimal
unstable foliations, see~\cite{ViY13} and also \cite{BoV00,BDP02,BDPP08};
accessible skew-products $M\times S^1 \to M\times S^1$ over Anosov diffeomorphisms which are
not rotation extensions, see~\cite{ViY13}. New examples will be given in Section~\ref{s.geo-explode}.


In what follows we give the precise statements of our results.

\subsection{Partial hyperbolicity, physical measures and skeletons}

In this paper, a diffeomorphism $f:M \to M$ is called \emph{partially hyperbolic}
if there is a continuous invariant splitting $TM=E^{cs} \oplus E^u$ of the
tangent bundle and there are constants $c>0$ and $\sigma>1$ such that
\begin{itemize}
\item[(a)] $\|Df^n v^u\| \ge c\sigma^n\|v^u\|$ for every $v^u\in E^u$ and every $n\ge 1$
(we say that $E^u$ is uniformly expanding).
\item[(b)] $E^{cs}$ is dominated by $E^u$:
$$
\frac{\|Df^n v^u\|}{\|Df^n v^{cs}\|} \ge c \sigma^n\frac{\|v^u\|}{\|v^{cs}\|}
$$
for every nonzero $v^u\in E^u$, $v^{cs}\in E^{cs}$, and every $n\ge 1$.
\end{itemize}
The \emph{unstable bundle} $E^u$ is automatically uniquely integrable:
there exists a unique foliation $\cF^u$ of $M$ with $C^1$ leaves tangent
to $E^u$ at every point. This \emph{unstable foliation} $\cF^u$ is invariant,
meaning that $f(\cF^u(x)) = \cF^u(f(x))$ for every $x\in M$ and the leaves
are, actually, as smooth as the diffeomorphism itself.


We call \emph{$u$-disk} any embedded disk contained in a leaf of the unstable foliation.
A partially hyperbolic map $f:M\to M$ has \emph{mostly contracting center}
(Bonatti, Viana~\cite{BoV00}) if, given any $u$-disk $D^u$, one has
$$
\limsup_{n\to\infty} \frac 1n \log \|Df^n(x) \mid E^{cs}(x)\|<0
$$
for every $x$ in some positive Lebesgue measure subset $D_0^u\subset D^u$.

A \emph{physical measure} for $f:M\to M$ is an invariant probability $\mu$
whose \emph{basin}
$$
B(\mu) = \{x\in M : \frac 1n \sum_{j=0}^{n-1} \delta_{f^j(x)} \text{ converges to $\mu$
in the  weak$^*$ topology}\}
$$
has positive volume. Bonatti, Viana~\cite{BoV00} proved that every $C^{1+\vep}$
diffeomorphism with mostly contracting center has a finite number of physical measures,
and the union of their basins contains almost every point in the ambient space.
See \cite{BDP02,Dol00} for several related results.
The set of Lebesgue density points of $B(\mu)$ will be called \emph{essential
basin} of $\mu$ and will be denoted $B_{ess}(\mu)$.

Let $f:M\to M$ be a $C^{1+\vep}$ partially hyperbolic diffeomorphism with mostly contracting center.
We say that a hyperbolic saddle point has \emph{maximum index} if the dimension of its
stable manifold coincides with the dimension of the center-stable bundle $E^{cs}$.
A \emph{skeleton} of $f$ is a collection $\cS=\{p_1,\cdots, p_k\}$ of hyperbolic
saddle points with maximum index satisfying
\begin{itemize}
\item[(i)] For any $x\in M$ there is $p_i \in \cS$ such that the stable manifold
$W^s(\Orb(p_i))$ has some point of transversal intersection with the unstable
leaf $\cF^u(x)$ through $x$.
\item[(ii)] $W^s(\Orb(p_i))\cap W^u(\Orb(p_j))=\emptyset$ for every $i\neq j$,
that is, the points in $\cS$ have no heteroclinic intersections.
\end{itemize}
Observe that a skeleton may not exist (for instance if $f$ has no periodic points).
Also, the skeleton needs not be unique, when it exists.
On the other hand, existence of a skeleton is a $C^1$-robust property, as
we will see in a while.

\begin{main}\label{t.mainA}
Let $f$ be a $C^{1+\vep}$ diffeomorphism with mostly contracting center. Then $f$ admits
some skeleton. Moreover, if $\cS=\{p_1, \ldots, p_k\}$ is a skeleton then for each
$p_i\in\cS$ there exists a distinct physical measure $\mu_i$ such that
\begin{enumerate}
\item the closure of $W^u(\Orb(p_i))$ and the homoclinic class of the orbit
$\Orb(p_i)$ both coincide with $\supp \mu_i$, which is the finite union of
disjoint $u$-minimal component, i.e., each unstable leaf in every component
is dense in this setting.
\item the closure of $W^s(\Orb(p_i))$ coincides with the closure of
the essential basin of the measure $\mu_i$.
\end{enumerate}
In particular, the number of physical measures is precisely $k=\#\cS$.
Moreover, $\supp(\mu_i)\cap \supp(\mu_j)=\emptyset$ for $1\leq i \neq j \leq k$.
\end{main}

In the proof (Section~\ref{s.srbgeostructure}) we just pick, for each physical measure
$\mu_i$ a hyperbolic periodic point $p_i\in\supp \mu_i$ with maximum index:
such points constitute a skeleton. When their stable manifolds are everywhere dense,
we get from part (b) of the theorem that there exist several physical measures,
whose basins are intermingled. Such examples, that generalize the main observation
of Kan~\cite{Ka94}, are exhibited in Section~\ref{s.geo-explode}.

\subsection{Variation of physical measures}

Theorem~\ref{t.mainA} provides us with a tool to mirror physical measures into
hyperbolic periodic points, and this can be used to describe the way physical
measures vary when the dynamics is modified.
Starting from a skeleton $\cS=\{p_1, \dots, p_k\}$ for $f$, we may consider its
continuation $\tilde\cS = \{p_1(g), \dots, p_k(g)\}$ for any nearby $g$.
Then any maximal subset of $\tilde\cS$ satisfying condition (ii) is a skeleton
for $g$. That is the main content of the following theorem:

\begin{main}\label{t.mainB}
There exists a $C^{1+\vep}$ neighborhood $\cU$ of $f$ such that, for any $g\in\cU$,
any maximal subset of the continuation $\{p_1(g), \ldots, p_k(g)\}$ which has
no heteroclinic intersections is a skeleton. Consequently, the number of physical
measures of $g$ is not larger than the number of physical measures of $f$.

In fact, these two numbers coincide if and only if there are no heteroclinic
intersections between the continuations $p_i(g)$. Moreover, in that case,
each physical measure of $g$ is close to some physical measure of $f$, in the
weak$^*$ topology.

In addition, restricted to any subset of $\cU$ where the number of physical
measures is constant, the supports of the physical measures and the closures
of their essential basins vary in a lower semi-continuous fashion with the
dynamics, both in the sense of the Hausdorff topology.
\end{main}

Of course, this implies that the number of physical measures is an upper semi-continuous
function of the dynamics. Consequently, this number is locally constant on an open and
dense subset of diffeomorphisms with mostly contracting center. These facts had been proved
before by Andersson~\cite{An10}. One important point in our approach is that we give a
definite explanation for possible ``collapse'' of physical measures: one physical measure
is lost for each  heteroclinic intersection that is created between the continuations of
elements of the skeleton. The precise statements are in Propositions~\ref{p.weakcollapse}
and~\ref{p.strongcollapse}.

We also want to explain how the basins of the physical measures vary with the
dynamics in the following measure theoretical sense. Define the pseudo-distance
$d(A,B)=\Leb(A \Delta B)$ in the space of measurable subsets of $M$.

\begin{main}\label{t.mainC}
Let $\cO$ be any subset of $C^{1+\vep}$ diffeomorphisms with mostly contracting center
such that all the diffeomorphisms in $\cO$ have the same number of physical measures.
Then their basins $B_i(f)$ vary continuously with $f\in\cO$, relative to the pseudo-distance $d$.
\end{main}

In Subsection~\ref{ss.examples} we will show how this theory can be applied to various examples,
including those of Kan~\cite{Ka94}. In particular, Theorem~\ref{t.mainC} shows that the
basins are quite stable from a measure-theoretical point of view.

\section{Geometric structure of physical measures \label{s.srbgeostructure}}

Let $f$ be a $C^{1+\vep}$ partially hyperbolic diffeomorphism with mostly contracting center.
As before, $E^{cs} \oplus E^{u}$ denotes the corresponding invariant splitting and $i_{cs}=\dim E^{cs}$.
We call \emph{Gibbs $u$-state} of $f$ any invariant probability absolutely continuous along
strong unstable leaves. It follows that the support is $u$-saturated, that is, it consists of
entire unstable leaves.

The notion of Gibbs $u$-state goes back to Pesin, Sinai~\cite{PS82} and was used by
Bonatti, Viana~\cite{BoV00} to construct the physical measures of diffeomorphisms with mostly
contracting center. Indeed, they showed that such diffeomorphisms have finitely many ergodic
Gibbs $u$-states, and these are, precisely, the physical measures.

Gibbs $u$-states also provide an alternative definition of mostly contracting center:
$f$ has mostly contracting center if and only if all Lyapunov exponents along the bundle
$E^{cs}$ are negative for every ergodic Gibbs $u$-state. This is related to the fact that,
given any disk $D$ inside an unstable leaf, any Cesaro accumulation point of the iterates of
(normalized) Lebesgue measure on $D$ is a Gibbs $u$-state. In fact, more is true:
every accumulation point of
$$
\frac 1n \sum_{j=0}^{n-1} \delta_{f^j(x)}
$$
is a Gibbs $u$-state, for almost every $x\in D$. Another useful property is that the space
$G(f)$ of all Gibbs $u$-states is convex and weak$^*$ compact. The extremal elements are the
ergodic Gibbs $u$-states. Moreover, $G(f)$ is an upper semi-continuous function of $f$,
in the sense that the set $\{(f,\mu): \mu \in G(f)\}$ is closed.
Proofs of these facts can be found in Chapter~11 of \cite{Beyond}.

%

The following fact will be used several times in what follows:

\begin{proposition}[Viana, Yang~\cite{ViY13}]\label{p.diffminimalfoliationinsupport}
If $f$ is a $C^{1+\vep}$ diffeomorphism with mostly contracting center then the supports
of its physical measures, $\mu_1, \dots, \mu_l$ are pairwise disjoint.
Moreover, the support of every $\mu_i$ has finitely many connected components and each
connected component is minimal for the unstable foliation (every unstable leaf is dense).
\end{proposition}

\subsection{Proof of Theorem~\ref{t.mainA}}

The first step is to construct a skeleton:

\begin{proposition}\label{p.existenceofgraph}
Every $C^{1+\vep}$ partially hyperbolic diffeomorphism with mostly contracting
center admits some skeleton.
\end{proposition}

\begin{proof}
Since the center Lyapunov exponents are all negative, every physical measure $\mu_i$, $1 \le i \le l$
is a hyperbolic measure (meaning that all the Lyapunov exponents are different from zero).
So, by Katok~\cite{Ka80}, there exist periodic points $q_i$ with maximum index and whose
stable manifold intersects transversely the unstable leaf of some point in the support of $\mu_i$.
Since the support is $u$-saturated, invariant, and closed, it follows that $q_i\in\supp\mu_i$.
For each $i$ we choose one such periodic point $q_i$; we are going to show that
$\{q_1, \dots, q_l\}$ is a skeleton for $f$.

Consider any $x\in M$ and let $D$ be a disk around $x$ inside the corresponding unstable leaf.
Let $\mu$ be any Cesaro accumulation point of the iterates of the volume measure $\Leb_D$ on $D$.
As observed before, $\mu$ is a Gibbs $u$-state and, hence, may be written as $\mu=\sum_{i=1}^l a_i \mu_i$.
Choose $i$ such that $a_i$ is non-zero. Let $B$ be a neighborhood of $q_i$ small enough that
the unstable leaf through any point in $B$ intersects the stable manifold $W^s(q_i)$ transversely.
Then $\mu_i(B)>0$, because $q_i\in\supp \mu_i$, and so $\mu(B)>0$.
Consequently, there is $n$ arbitrarily large such that $(f_*^n \Leb_D)(B)>0$.
This implies that the unstable manifold of $f^n(x)$ intersects $W^s(q_i)$ transversely.
By invariance, it follows that $\cF^u(x)$ intersects transversely the stable manifold of some iterate of $q_i$.
This proves condition (i) in the definition of skeleton.

Condition (ii) is easy to prove. Indeed, on the one hand, $W^u(\Orb(q_i))$ is contained in $\supp\mu_i$.
On the other hand, this support can not intersect $W^s(\Orb(q_j))$ for any $j\neq i$: otherwise,
$q_j$ would be in $\supp\mu_i$, which would contradict the fact that the supports are pairwise disjoint.
Thus, there can indeed be no heteroclinic connections.
\end{proof}

Now, we use the skeleton to analyse the physical measures:

\begin{proposition}\label{p.numberofmeasures}
Let $f$ be a $C^{1+\vep}$ diffeomorphism with mostly contracting center.
Suppose that $\cS=\{p_1,\cdots,p_k\}$ is a skeleton of $f$. Then
\begin{itemize}
\item[(a)] $\# \cS$ coincides with the number of physical measure of $f$;

\item[(b)] the closure of $W^u(Orb(p_i))$ coincides with $\supp(\mu_i)=H(p_i,f)$;

\item[(c)] the closure of $W^s(Orb(p_i))$ coincides with the closure of $B_{ess}(\mu_i)$.
\end{itemize}
\end{proposition}

\begin{proof}
To prove claim (a) it suffices to show that all skeletons have the same number of elements
(because the claim holds for the skeleton constructed in Proposition~\ref{p.existenceofgraph}).
Let $\cS'=\{q_1, \dots, q_l\}$ be any other skeleton.
By condition (i) in the definition, for each $q_j\in\cS'$ there is some $p_i\in\cS$
such that $W^u(q_j)$ intersects $W^s(\Orb(p_i))$ transversely.
Choose any such $p_i$ (we will see in a while that the choice is unique).
For the same reason, for this $p_i$ there exists some $q_t\in \cS'$ such that $W^u(p_i)$
intersects $W^s(\Orb(q_t))$ transversely.
It follows that $W^u(q_j)$ accumulates on $\Orb(q_t)$ which, by condition (ii) in the
definition, can only happen if $q_j=q_t$.
Thus, $p_i$ and $q_j$ are heteroclinically related to one another.
Since different elements of either skeleton do not have heteroclinic intersections,
this implies that $p_i$ is unique and the map $q_j \mapsto p_i$ is injective.
Reversing the roles of the two skeletons, we also get an injective map $p_i \mapsto q_j$
which, by construction, is the inverse of the previous one.
Thus, these maps are bijections and, in particular, $\#\cS=\#\cS'$.

Now take $\cS'$ to be the skeleton obtained in Proposition~\ref{p.existenceofgraph}.
Up to renumbering, we may assume that the $i=j$ in the previous construction.
Also by construction, each $p_i$ is contained in the closure of $W^u(\Orb(q_i))$,
which coincides with the support of $\mu_i$. Since the unstable foliation is minimal
in each connected component of the support, this implies that the closure of
$W^u(\Orb(p_i))$ coincides with $\supp(\mu_i)$. To finish the proof of claim (b)
it remains to show that this coincides with the homoclinic class of $p_i$.
We only have to prove that $H(p_i)$ contains the closure of $W^u(\Orb(p_i))$,
since the converse is an immediate consequence of the definition of homoclinic class.

To this end, let $D$ be any disk contained in the unstable manifold of $\Orb(p_i)$.
Let $\mu$ be any Cesaro accumulation point of the iterates $f^n_*\Leb_D$.
This is a Gibbs $u$-state and it gives full measure to $\supp \mu_i$ (because
$W^u(\Orb(p_i)) \subset \supp\mu_i$ and the latter is a compact invariant set).
Given that there are finitely many ergodic Gibbs $u$-states, and their supports are disjoint,
this implies that $\mu = \mu_i$. Then, by the same argument that we used in the
previous proposition, there exists some large $m$ large such that $f^m(D)$ intersects
$W^s(p_i)$ transversely. Since $D$ is arbitrary, this means that homoclinic points
are dense in the unstable manifold of $\Orb(p_i)$, which implies the claim.

It remains to prove the claim (c). Let $D$ be any disk contained in $W^u(p_i)$.
By Theorem~11.16 in \cite{Beyond}, for Lebesgue almost every $x\in D$ the sequence
$$
\frac 1n \sum_{j=0}^{n-1} \delta_{f^j(x)}
$$
converges to some Gibbs $u$-state. This Gibbs $u$-state must be $\mu_i$, because
by Proposition~\ref{p.diffminimalfoliationinsupport} and part (b) of
Proposition~\ref{p.numberofmeasures}, this is the unique ergodic Gibbs $u$-state that gives
weight to the closure of $W^u(\cO(p_i))$.
This proves that the basin of $\mu_i$ intersects $D$ on a full Lebesgue measure subset.

We claim that there exists a positive Lebesgue measure subset $D_1$ inside that intersection
such that the stable set of any point $y\in D_1$ contain an $i_{cs}$-dimensional disk
$W^s_{loc}(y)$ with uniform size; moreover, these local stable disks constitute an
absolutely continuous lamination (that is, the holonomy maps of this lamination preserve
zero measure sets). Indeed, let $\Lambda$ be any compact (non-invariant) set
with $\mu(\Lambda)>0$ such that every point in $\Lambda$ has a Pesin stable manifold
with uniform size, and these stable manifolds constitute an absolutely continuous
lamination (existence of such sets is a classical fact in Pesin theory~\cite{Pes76}).
It follows from the previous paragraph that the forward trajectory of almost every $x\in D$
accumulates on $\Lambda$. Thus one can find a neighborhood $V$ of $x$ inside $D$ such that
some large iterate $f^n(V)$ intersects $\Lambda$ on a positive Lebesgue measure subset.
Just take $D_1=V \cap f^{-n}(\Lambda)$. See also \cite[Lemma~6.6]{ViY13} for a similar
statement.

Let $x_0\in D_0$ be a Lebesgue density point for $D_1$ inside $D_0$.
Since the basin contains the stable sets of all points in $D_1$, and these are transverse to $D_0$,
it follows that every point in the local stable disk of $x_0$
is also a Lebesgue density for the basin in
ambient space. In particular, $W^s_{loc}(x_0)$ is contained in $B_{ess}(\mu_i)$.
Since $f^{-n}(W^s_{loc}(x_0)$ accumulates on $W^s(\Orb(p_i))$ and the essential basin is $f$-invariant,
it follows that $W^s(\Orb(p_i))$ is contained in the closure of $B_{ess}(\mu_i)$.

Now we prove the converse inequality. Let $x_0$ be any Lebesgue density point of the basin of $\mu_i$
in ambient space. Using the fact that the unstable foliation is absolutely continuous (see~\cite{BP74}),
we can find a small disk $D$ around $x_0$ inside the corresponding unstable leaf such that
$Leb_D(D\cap B(\mu_i))>0$. Let $B$ be a neighborhood of $p_i$ small enough that $\cF^u(y)$ intersects
$W^s(p_i)$ transversely, for every $y\in B$. Take $x \in D \cap B(\mu_i)$.
While proving part (b) we have shown that for such a point there exists arbitrarily large values of $n\ge 1$
such that $f^n(x) \in B$. Then $f^n(D)$ intersects $W^s(p_i)$ transversely and, hence, $W^s(\Orb(p_i))$
intersects $D$. Since $D$ is arbitrary, it follows that $x_0$ is in the closure of $W^s(\Orb(p_i))$.
\end{proof}

Combining Propositions~\ref{p.existenceofgraph} and~\ref{p.numberofmeasures} yields Theorem~\ref{t.mainA}.

\subsection{Proof of Theorem~\ref{t.mainB}}

It will be convenient to separate the two conditions in the definition of skeleton.
Let us call \emph{pre-skeleton} any finite collection $\{p_1,\cdots, p_k\}$ of saddles with maximum index
satisfying condition (i), that is, such that every unstable leaf $\cF^u(x)$ has some point of
transverse intersection with $W^s(\Orb(p_i))$ for some $i$. Thus a pre-skeleton is a skeleton if and only
if there are no heteroclinic intersections between any of its points.

One reason why this notion is useful is that the continuation of a pre-skeleton is always a pre-skeleton:

\begin{lemma}\label{l.robustgraph}
Let $f$ be a partially hyperbolic diffeomorphism which has a pre-skeleton $\cS=\{p_1,\cdots, p_k\}$.
Let $p_i(g)$, $i=1, \dots, k$ be the continuation of the saddles $p_i$ for nearby diffeomorphism $g$.
Then $\cS(g)=\{p_1(g),\cdots, p_k(g)\}$ is a pre-skeleton for every $g$ in a neighborhood of $f$.
\end{lemma}

\begin{proof}
This is a really a simple consequence of the fact that the unstable foliation depends continuously on the
point and the dynamics. Let us detail the argument.
Given any $x\in M$, take $i$ such that the unstable leaf $\cF^u(x)$ has some transverse intersection $a_x$
with the stable manifold of some point in the orbit of $p_i\in\cS$.
Fix $R_x>0$ large enough so that $a_x$ is in the interior of the $R_x$-neighborhood $\cF_{R_x}^u(x)$ of $x$
inside $\cF^u(x)$ and in the interior of the $R_x$-neighborhood $W^s_{R_x}(\Orb(p_i))$ of the orbit of $p_i$
inside its stable manifold. Then, since unstable leaves vary continuously with the point,
for any $y$ in a small neighborhood $U_x$ of $x$, there exists $a_y$ close to $a_x$ such that $\cF^u_{R_x}(y)$
and $W_{R_x}^s(\Orb(p_i))$ intersect transversely at $a_y$.
Let $\{U(x_1), \dots, U(x_m)\}$ be a finite covering of $M$ and let $R=\max\{R_{x_1}, \dots, R_{x_m}\}$.
Thus, $\cF_R^u(x)$ has some transverse intersection with $\cup_{i=1}^k W^s_R(\Orb(p_i))$ for every $x\in M$.
Since unstable leaves also vary continuously with the dynamics, it follows that there is a $C^1$ neighborhood
$\cU$ of $f$ such that $\cF^u_R(x,g)$ has some transverse intersection with $\cup_{i=1}^k W^s_R(p_i(g))$
for every $x\in M$ and every $g\in\cU$.
\end{proof}

Another reason why the notion of pre-skeleton is useful to us is that every pre-skeleton contains some
skeleton. To prove this it is convenient to introduce the following partial order relation, which will also be
useful later on. For any two elements of a pre-skeleton $\cS=\{p_1, \dots, p_k\}$ define:
$p_i \prec p_j$ if and only if
$$
W^u(\Orb(p_i)) \text{ has some transverse intersection with } W^s(\Orb(p_j)).
$$
It follows from the inclination lemma of Palis (\cite[\S~7]{PM82}) that $\prec$ is transitive and thus a
partial order relation. We say that $p_i\in\cS$ is a \emph{maximal element} if $p_j \prec p_i$ for every
$p_j\in\cS$ such that $p_i\prec p_j$.
Two maximal elements $p_i$ and $p_j$ are \emph{equivalent} if $p_i \prec p_j$ and $p_j \prec p_i$.
We call \emph{slice} of $\cS$ any subset that contains exactly one element in each equivalence class
of maximal elements.

\begin{lemma}\label{l.existenceskeleton}
Let $f$ be a partially hyperbolic diffeomorphism which has a pre-skeleton $\cS=\{p_1,\cdots, p_k\}$.
Any slice of $\cS$ is a skeleton.
\end{lemma}

\begin{proof}
Let $\cS'$ be a subset as in the statement. Begin by noting that $\cS'$ is also a pre-skeleton.
Indeed, since $\cS$ is assumed to be a pre-skeleton, for any $x\in M$ there exists $p_i\in\cS$ such that
$\cF^u(x)$ has some transverse intersection with $W^s(\Orb(p_i))$. Moreover, there exists some maximal
element $p_j$ of $\cS$ such that $p_i \prec p_j$. Using the $\lambda$-lemma, it follows that
$\cF^u(x)$ has some transverse intersection with $W^s(\Orb(p_j))$. Moreover, up to replacing $p_j$ by
some other maximal element equivalent to it, we may suppose that $p_j \in \cS'$. This proves our claim.
Finally, by definition, there is no heteroclinic intersection between the elements of $\cS'$.
So, $\cS'$ is indeed a skeleton.
\end{proof}

Now we are ready to give the proof of Theorem~\ref{t.mainB}. The set $\cS=\{p_1,\cdots,p_k\}$ is a
pre-skeleton of $f$, of course. So, by Lemma~\ref{l.robustgraph}, there is a $C^1$ neighborhood
$\cV$ of $f$ such that $\cS(g)=\{p_1(g), \dots, p_k(g)\}$ is a pre-skeleton for every $g\in\cV$.
Since diffeomorphisms with mostly contracting center form a $C^{1+\vep}$ open set (by Andersson~\cite{An10}),
we may find a $C^{1+\vep}$ neighborhood $\cU \subset\cV$ such that every $g\in\cU$ has mostly
contracting center. By Lemma~\ref{l.existenceskeleton}, every slice $\cS'(g)$ of $\cS(g)$ is a
skeleton for $g$. Since $\#\cS'(g) \le \# \cS(g) = \cS(f)$, it follows from Theorem~\ref{t.mainA}
that  the number of physical measures of $g\in\cU$ is not larger than the number of physical
measures of $f$. Indeed, these two numbers coincide if and only if $\cS(g) $ is a skeleton for $g$,
that is, if there are no heteroclinic intersections between the continuations $p_i(g)$.
This proves the first part of the theorem.

Now let $(f_n)_n$ be a sequence of diffeomorphisms converging to $f$ in the $C^{1+\vep}$ topology
and suppose that $\cS(f_n)=\{p_1(f_n),\cdots, p_k(f_n)\}$ is a skeleton of $f_n$ for any large $n$.
Let $\mu_1(f_n), \dots, \mu_k(f_n)$ be the physical measures (ergodic Gibbs $u$-states).
By Theorem~\ref{t.mainA}, we may number these measures in such a way that each $\mu_i$ is supported
on the closure of $W^u(\Orb(p_i(f_n)))$. Up to restricting to a subsequence, we may assume that
$\mu_i(f_n)$ converges, in the weak$^*$ topology, to some $f$-invariant measure $\mu_i^*$.
By semicontinuity of the space of Gibbs $u$-states, every $\mu_i^*$ is a Gibbs $u$-state for $f$.
Write $\mu$ as a convex combination $\mu_i^*=\sum_{j=1}^k a_j\mu_j$ of the physical measures of $f$.
We claim that $a_i=1$. Indeed, suppose that there is $j \neq i$ such that $a_j>0$. Then
$$
\limsup_n \supp(\mu_i(f_n)) \supset \supp(\mu_j(f)).
$$
By Theorem~\ref{t.mainA}, we have that $\supp(\mu_i(f_n))=$ closure of $W^u(Orb(p_i),f_n)$.
For $n$ large, this implies that $W^u(\Orb(p_i(f_n)),f_n)$ has some transverse intersection with
$W^s_{loc}(\Orb(p_j(f)),f)$, because the unstable manifolds of hyperbolic periodic points vary
continuously with the dynamics. Using the corresponding fact for stable manifolds, we
conclude that  $W^u(\Orb(p_i(f_n)),f_n)$ has some transverse intersection with
$W^s_{loc}(\Orb(p_j(f_n)),f_n)$. This contradicts the assumption that $S(f_n)=\{p_1(f_n),\cdots,p_k(f_n)\}$
is a skeleton of $f_n$. This proves our claim, which yields the second part of the theorem.

By the stable manifold theorem (see \cite[Theorem~6.2]{PM82} and \cite[page~154]{PT93}),
for each $R>0$,
the local invariant manifolds $W^s_R(\Orb(p_i(g)))$ and $W^u_R(\Orb(p_i(g)))$ vary continuously with $g.$
This implies that their closures vary in a lower semi-continuous fashion with $g$,
relative to the Hausdorff topology. By parts (b) and (c) of Proposition~\ref{p.numberofmeasures},
this means that both the supports and the closures of the essential basins of the physical measures
vary lower semi-continuously with the dynamics, as claimed in the third part of the theorem.

The proof of Theorem~\ref{t.mainB} is complete.

\subsection{Local description of the continuation of physical measures}

Our next goal will be to analyse how physical measures and their basins vary with the dynamics.
Here we find a couple of conditions that ensure continuous dependence. This is a prelude to the
next section, where we will analyse how physical measures may collapse as their basins explode.

Take $f$ to be a diffeomorphism with mostly contracting center with a skeleton $\cS=\{p_1, \dots, p_k\}$.
Let $\cS(g)=\{p_1(g), \dots, p_k(g)\}$ be its continuation for nearby diffeomorphisms $g$.

\begin{corollary}\label{c.localsinglestable}
Let $i\in\{1, \dots, k\}$ and $(f_n)_n$ be a sequence converging to $f$ in $\Diff^{1+\vep}(M)$
such that for every $n$ the point $p_i(f_n)$ is a maximal element of $\cS(f_n)$ and no other
element of $\cS(f_n)$ is equivalent to $p_i(f_n)$. Then each $f_n$ has a physical measure
$\mu_i(f_n)$ on the closure of $W^u(Orb(p_i(f_n)))$ such that these physical measures converge
to $\mu_i$ in the weak$^*$ topology as $n\to\infty$.
\end{corollary}

\begin{proof}
By Lemma~\ref{l.existenceskeleton}, each $f_n$ admits a physical measure $\mu_i(f_n)$ supported
on the closure of the unstable manifold of $\Orb(p_i(f_n))$.
Suppose that $(\mu_i(f_n))_n$ does not converge to $\mu_i$.
We may assume that the sequence converges to some measure $\mu$. Then $\mu$ is a Gibbs $u$-state
of $f$ and so we may write it as
$$
\mu=a_1\mu_1 + \dots + a_i\mu_i + \dots + a_k \mu_k.
$$
Since $\mu \neq \mu_i$, there exists $j\neq i$ such that $a_j \neq 0$.
By the same argument as in the proof of Theorem~\ref{t.mainB}, we have that $W^u(\Orb(p_i(f_n)),f_n)$
intersects $W^s(\Orb(p_j),f)$ transversely at some point, for every large $n$.
Consequently, if $n$ is large enough then $W^u(\Orb(p_i(f_n)),f_n)$ has some transverse intersection
with $W^s(\Orb(p_j(f_n)), f_n)$. This implies that $p_i(f_n) \prec p_j(f_n)$, which contradicts the
assumption that $p_i(f_n)$ is maximal and its equivalence class is formed by a single point.
\end{proof}

Given $r\geq 1$ and two saddle points $p(f)$ and $q(f)$ of diffeomorphism $f$, we say that $q$ is not \emph{$C^r$ attainable}
from $p$ if there is a $C^r$ neighborhood $\cV$ of $f$ such that $W^u(p(g),g)\cap W^s(q(g),g)=\emptyset$
for any $g\in \cV$, where $p(g)$ and $q(g)$ are the analytic continuations of $p(f)$ and $q(f)$,
respectively.

\begin{corollary}\label{c.localrobuststable}
Assume that $p_i(f)\in\cS$ is not $C^{1+\vep}$ attainable from any $p_j(f)\in\cS$ with $j\neq i$.
Then the physical measure $\mu_i(f)$ is stable, in the sense that for every $g$ in a $C^{1+\vep}$ neighborhood of $f$
there exists a physical measure $\mu_i(g)$ which is close to $\mu_i(f)$ in the weak$^*$ topology.
\end{corollary}

\begin{proof}
Let $\cV$ be a neighborhood of $f$ as in the definition of non-attainability. Let $\cS'(g)$ be any
slice of $\cS(g)$. By Lemma~\ref{l.existenceskeleton}, $\cS'(g)$ is a skeleton for $g$.
The assumption implies that $p_i(g)$ is a maximal element of $\cS'(g)$ and its equivalence class
consists of a single point. So, the conclusion follows from Corollary~\ref{c.localsinglestable}.
\end{proof}

\section{Exploding basins\label{s.geo-explode}}

We start by giving a geometric and measure-theoretical criterion for a partially hyperbolic diffeomorphism
to have mostly contracting center, using the notion of skeleton and a local version of the mostly contracting
center property. Then we use this criterion to give new examples of diffeomorphisms with any
finite number of physical measures, whose basins are all dense in the ambient space.

Such examples are not stable: the number of physical measures may decrease under perturbation.
Indeed, for any proper subset of physical measures one can find a small perturbation of the original
diffeomorphism for which those physical measures disappear
(their basins are engulfed by the basins of the physical measures that do remain).

Using different perturbations, one can approximate the original diffeomorphism $f$ by other diffeomorphisms
$f_n$ having a unique physical measure $\mu_n$, in such a way that $(\mu_n)_n$ converges to any given
Gibbs $u$-state of $f$. In particular, such examples are \emph{statistically unstable}: the simplex
generated by all the physical measures does not vary continuously.

\subsection{Criterion}

Take $f$ to be a partially hyperbolic diffeomorphism with invariant splitting $E^u \oplus E^{cs}$.
As before, denote $i_{cs}=\dim E^{cs}$.
We start with a semi-local version of the notion of mostly contracting center.

Let $\Lambda$ be a compact $u$-saturated $f$-invariant subset of $M$.
We say that \emph{$f$ has mostly contracting center at $\Lambda$} if the center Lyapunov exponents
are negative for every ergodic Gibbs $u$-state supported on $\Lambda$.
Then, we say that $\Lambda$ is \emph{an elementary set} if there exists exactly one ergodic Gibbs $u$-state
$\mu$ supported in $\Lambda$ and it satisfies $\supp\mu=\Lambda$.

The same arguments as in Theorem~\ref{t.mainA} also yield a corresponding semi-local statement:
If $\Lambda$ is an elementary set and $\mu$ is the corresponding Gibbs $u$-state, then
\begin{itemize}
\item $\mu$ is a physical measure;
\item $\Lambda$ has finitely many connected components and the unstable foliation is minimal in each connected component;
\item if $p\in \Lambda$ is any hyperbolic saddle with maximum index, then the closure
of $W^s(\Orb(p))$ coincides with the closure of the essential basin of $\mu$.
\end{itemize}
$\Lambda$ contains some hyperbolic saddle with maximum index, by arguments in the proof of Proposition~\ref{p.existenceofgraph}.

\begin{proposition}\label{p.newcretetion}
Let $\Lambda_1$ , \dots, $\Lambda_k$ be pairwise disjoint elementary sets,
$\mu_1$, \dots, $\mu_k$ be the corresponding Gibbs $u$-states,
and $p_i\in\Lambda_i$, $i=1, \dots, k$ be hyperbolic saddles with maximum index.
If $\{p_1, \dots, p_k\}$ is a pre-skeleton, then it is a skeleton, and $f$ has mostly contracting center.
Moreover, $\{\mu_1, \dots, \mu_k\}$ are the physical measures of $f$, and their basins cover a full Lebesgue measure subset.
\end{proposition}

\begin{proof}
If some unstable manifold $W^u(\Orb(p_i))$ intersects some stable manifold $W^s(\Orb(p_j))$ then,
by the inclination lemma~\cite[\S~7]{PM82}, the closure of $W^u(\Orb(p_i))$ intersects the closure
of $W^u(\Orb(p_j))$. By the definition of elementary sets, this implies that $\Lambda_i$ intersects
$\Lambda_j$ and, in view of our assumptions, that can only happen if $i=j$. This proves that
$\{p_1, \dots, p_k\}$ is a skeleton.

Now let us check that $f$ has mostly contracting center.
It is part of the definition of elementary set that the center Lyapunov exponents of $\mu_j$
are all negative, for every $j=1, \dots, k$. So, to prove that $f$ has mostly contracting
center it suffices to show that $f$ has no any other ergodic Gibbs-$u$ states.
Suppose there exists some ergodic Gibbs-$u$ state $\mu\notin \{\mu_1,\cdots,\mu_k\}$.
It follows from the definition that there exists a $u$-disk $D$ contained in some unstable
leaf that intersects the basin of $\mu$ on a full Lebesgue measure set $D_0\subset D$.
We claim that there exist $n_0\ge 1$ and $1\leq i \leq k$ such that $f^{n_0}(D_0)$ intersects
the basin of $\mu_i$. Of course, this contradicts the fact that $\mu \neq \mu_i$.
Thus, we are left to justify our claim.

Since $\{p_1,\cdots, p_k\}$ is a pre-skeleton, there exist $n \ge 1$ and $1\leq i \leq k$
such that $f^n(D)$ intersects $W^s(p_i)$ transversely at some point (otherwise, the
Hausdorff limit of $f^n(D)$ would contain some unstable leaf disjoint from
$\cup_{i=1}^k W^s(\cO(p_i))$, which would contradict the definition of pre-skeleton).
Again by the definition of Gibbs $u$-state, there exists a $u$-disk $D'\subset \supp\mu_i$
and a full Lebesgue measure subset $D'_0\subset D'$ formed by regular points of $\mu_i$.
Since the center Lyapunov exponents are negative, it follows from Pesin theory that there
exists a lamination whose laminae are local stable manifolds $W^s_{loc}(x)$ of almost
every point $x \in D_0'$. Moreover, this \emph{stable lamination} is absolutely continuous.

Theorem~11.16 in \cite{Beyond} gives that the time average of Lebesgue almost every
$x\in W^u(p_i)$ is a Gibss $u$-state. By the definition of  elementary set, this Gibbs
$u$-state must be $\mu_i$. Moreover, the orbit of any such $x$ must accumulate on the whole
$\supp\mu_i=\Lambda_i$. In particular, $W^u(\Orb(p_i))$ is dense in $\supp\mu_i$.
Assuming that $n_0$ is large enough, $f^{n_0}(D)$ is close to $W^u(\Orb(p_i))$ and,
in particular, it cuts $\cup_{x\in D'_0}W^s_{loc}(x)$. The intersection is contained in the
basin of $\mu_i$, since $W^s_{loc}(x) \subset B(\mu_i)$ for every $x \in D_0'$.
Moreover, by absolute continuity of the lamination, the intersection has positive
Lebesgue measure. This implies that $f^n(D_0)$ intersects the basin of $\mu_i$.
\end{proof}

\subsection{New Kan-type examples}\label{Kan's example with several measures}\label{Kan's example}

In this subsection, we use Proposition~\ref{p.newcretetion} to construct new examples of
diffeomorphisms with mostly contracting center and several physical measures, such that
every basin intersects every open set on a positive measure subset.

\begin{proposition}\label{p.Kanseveralmeasures}
For any $k\ge 1$, there is a diffeomorphism $f\in \Diff^2(T^2\times S^2)$ such that
$f$ has mostly contracting center and $k$ physical measures $\mu_1, \dots, \mu_k$
such that $\supp \mu_i = T^2\times A_i$ for some $A_i\subset S^2$ and the basin
$B(\mu_i)$ is dense in $T^2 \times S^2$, for every $i$. Moreover, the same remains
true for any diffeomorphism in a $C^2$-neighborhood which preserves the set
$T^2\times A_i$ for all $i=1, \dots k$.
\end{proposition}

\begin{proof}
Let $k$ be fixed and $g\in \Diff^1(T^2)$ be a $C^r$ Anosov diffeomorphism with $2k$ fixed points,
denoted as $p_1, p'_1, \dots, p_k, p'_k$. Our example will be a partially hyperbolic skew product map
$$
f: T^2 \times S^2 \rightarrow T^2\times S^2,
\quad f(x,y)=(g(x), h_x(y)),
$$
whose center foliation is the vertical foliation by spheres, $W^c(x)=\{x\}\times S^2$.
It is easy to see that, for any $x\in T^2$,
$$
W^s(W^c(x), f)=W^s(x,g) \times S^2
\quad\text{and}\quad
W^u(W^c(x), f)=W^u(x,g) \times S^2.
$$
For $x$ and $\tilde x$ in the same stable manifold of $g$,
let $H^s_{x,\tilde x}:W^c(x) \to W^c(\tilde x)$ be the \emph{stable holonomy},
defined as the projection along strong stable leaves of $f$.
Let the \emph{unstable holonomy} $H^u_{x,\tilde x}:W^c(x) \to W^c(\tilde x)$ be defined analogously,
for $x$ and $\tilde x$ in the same unstable leaf of $g$.

Assuming that $h_x$ is uniformly close to the identity in the $C^2$ topology,
the partially hyperbolic map $f$ is \emph{center bunched} (see~\cite{PSh97} or~\cite{BW10}),
so that these holonomy maps are all $C^1$ diffeomorphisms; moreover, they are close to the identity
in the $C^1$ topology. In what follows we consider $k\ge 3$: the cases $k=1, 2$ are easier.

For the time being, take $k$ to be even; the odd case will be treated at the end of the proof.
Let $\cC$ and $\cC'$ be two smooth closed curves in $S^2$ intersecting transversely on exactly $k$ points,
$A_1, \dots, A_k$.  Take these points to be listed in cyclic order.
Then consider points $B_1, \dots, B_k \in \cC$, such that each $B_i$ lies in the circle segment
between $A_i$ and $A_{i+1}$ (with $A_{k+1}=A_1$). For each $i=1, \dots, k$, let $X_i$ be a Morse-Smale
vector field on the sphere such that:
\begin{enumerate}
\item[(i)] $\Omega(X_i) = \{A_1, B_1, \dots, A_k, B_k\}$;
\item[(ii)] $A_i$ is a sink, $A_1, \dots, A_{i-1}, A_{i+1}, \dots, A_k$ are saddles
and $B_1, \dots, B_k$ are sources;
\item[(iii)] the basin of the attractor $A_i$ is the complement of segment $S_i\subset\cC$
connecting all the saddles and sources.
\end{enumerate}
Figure~\ref{f.flow} illustrates the case $k=4$ and $i=1$: then $S_1$ is just the segment
of $\cC$ from $B_1$ to $B_4$ that does contain $A_1$.

\begin{figure}[ht]
\begin{center}
\psfrag{C}{$\cC$}
\psfrag{A1}{$A_1$}\psfrag{A2}{$A_2$}\psfrag{A3}{$A_3$}\psfrag{A4}{$A_4$}
\psfrag{B1}{$B_1$}\psfrag{B2}{$B_2$}\psfrag{B3}{$B_3$}\psfrag{B4}{$B_4$}
\includegraphics[width=2in]{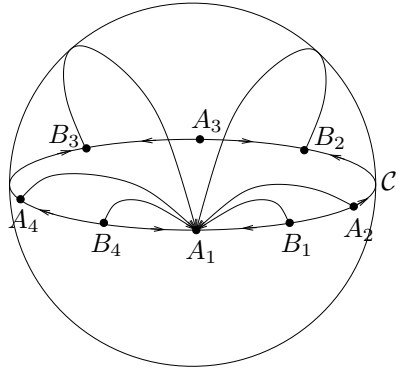}
\caption{\label{f.flow}Morse-Smale vector field on the sphere}
\end{center}
\end{figure}

Analogously, consider points $B'_1, \dots, B'_k \in \cC'$, such that each $B'_i$ lies in the segment
of $\cC'$ between $A_i$ and $A_{i+1}$. Then let $X'_i$, $i=1, \dots, k$, be a Morse-Smale
vector field on the sphere satisfying (i), (ii) and (iii), with $B_i$ replaced by $B_i'$ and $S_i$
replaced by a segment $S_i'\subset\cC'$. Let us consider a partially hyperbolic skew-product
$f:\TT^2\times S^2 \to \TT^2 \times S^2$ satisfying
\begin{itemize}
\item[(1)] $A_1, \dots, A_k$ are fixed points of $h_x(\cdot)$ for any $x\in T^2$.
\item[(2)] $h_{p_i}=$ time-$\vep$ map of $X_i$ and $h_{p'_i}=$ time-$\vep$ of $X_i'$,
for some small $\vep>0$.
\item[(3)] $f$ is $C^2$ close to $(g(x),id)$.
\end{itemize}
Condition (1) means that each $T_i=\TT^2\times\{A_i\}$, $i=1, \dots, k$ is an $F$-invariant
torus; clearly, the restriction $F \mid T_i$ is an Anosov map.
It is also clear that the three conditions are compatible, as long as we choose $\vep$ in (2)
sufficiently small. For example, we may take $h_x$ to be the identity map on $S^2$ for every $x$
outside small neighborhoods of $p_1$, \dots, $p_k$ and $p_1'$, \dots, $p_k'$.
Then, we may modify these maps $h_x$ to make them contracting at each $A_i$ (preserving the
previous three conditions), so that
\begin{itemize}
\item[(4)] $\int_{T_i} \log \|Dh_x(A_i)\| \, d\mu_i(x)<0$ for $i=1, \dots, k$,
where $\mu_i$ denotes the (unique) Gibbs $u$-state of $F \mid T_i$.
\end{itemize}
This last condition implies that the center Lyapunov exponents of every $\mu_i$ are negative,
and so $T_i$ is an elementary set.

\begin{lemma}\label{l.skeleton}
The set $\{p_i\times A_i\}_{i=1}^k$ is a skeleton.
\end{lemma}

\begin{proof}
As a first step, we prove that every strong unstable leaf $\cF^u(z)$ has a point of transverse intersection
with the stable manifold of some $(p_i,A_i)$.
Observe that $W^s(W^c(p_1),f)=W^s(p_1,g) \times S^2$. Also, $W^s(p_1,g)$ intersects the $g$-unstable manifold
of any point in $\TT^2$ transversely (recall that $g:\TT^2\to\TT^2$ is Anosov).
It follows that, for any $z\in M$, $\cF^u(z)$ intersects $W^s(W^c(p_1),f)$ transversely at some point $a$.
There are three possibilities:
\begin{itemize}
\item[(a)] $a\in W^s((p_1,A_1),f)$;
\item[(b)] $a\in W^s((p_1, A_i),f)$ for some $i \neq 1$;
\item[(c)] $a\in \cF^s(B_j)$ for some $1\leq j\leq k$.
\end{itemize}

In case (a) we are done. As for case (b), we claim that it implies that $\cF^u(z)$ has some transverse
intersection $W^s((p_i, A_i),f)$. Indeed, the hypothesis implies that the iterates $f^n(\cF^u(z))$
accumulate on the unstable leaf $\cF^u(p_1,A_i)$ of the fixed point $(p_1, A_i)$.
The latter is contained in the Anosov torus $T_i$, which also contains $(p_i,A_i)$ and its strong
stable leaf $\cF^s(p_i,A_i)$. In fact, $\cF^u(p_1,A_i)$ and $\cF^s(p_i,A_i)$ are transverse inside $\TT_i$.
Thus, it follows that the iterates $f^n(\cF^u(z))$ accumulate on $(p_i,A_i)$.
Since $(p_i,A_i)$ has stable index $3$, we get that $f^n(\cF^u(z))$  has some transverse intersection
with $W^s((A_i,p_i),f)$ for every large $n$. Taking pre-images, we get our claim.
Thus, in case (b) we are done as well.

Now, we consider case (c). For $n$ large, $f^n(\cF^u(z))=f^n(\cF^u(a))=\cF^u(f^n(a))$ is
close to $\cF^u(p_1,B_j)$. Let $q\in\TT^2$ be a point of transverse intersection between
$W^u(p_1, g)$ and $W^s(\tilde{p}_1,g)$.
Then the consider the map
$$
H = H^s_{q,\tilde{p}_1} \circ H^u_{p_1,q}: W^c(p_1) \rightarrow W^c(\tilde{p}_1).
$$
As observed above, under our assumptions the map $H$ is $C^1$ close to the identity map in
the second coordinate. So, in view of our conditions on $\cC$ and $\cC'$ (more specifically,
the assumption that they meet at $A_1, \dots, A_k$ only, and they do so transversely), we have
that $H(B_j) \notin \cC'$. Consequently, $H(B_j)\in W^s((\tilde{p}_1,A_1)$. This means that
the strong unstable leaf $\cF^u(p_1,B_j)$ has some transverse intersection with
$W^s((\tilde{p}_1,A_1),f)$. Then the same is true for $f^n(\cF^u(z))$ if $n$ is large enough.
Now observe that $(p_1,A_1)$ and $(\tilde{p}_1,A_1)$ are homoclinically related, meaning that
the unstable manifold of any point has some transverse intersection with the stable manifold
of the other. So, the previous conclusion implies that $f^n(\cF^u(z))$ has some transverse
intersection with $W^s((p_1,A_1),f)$. This reduces the present situation to case (a).

Thus, we have shown that $\{p_i\times A_i\}_{i=1}^k$ is a pre-skeleton. Next, notice that
$W^u((p_i,A_i),f)=\cF^u(p_i\times A_i)$ is contained in $T_i$ for every $i$. Since these tori
are pairwise disjoint, and each one of them is fixed under $f$, we have that $\TT_i$ is in
the complement of $W^s((p_j,A_j),f)$ for every $j \neq i$. So, the points $(p_i,A_i)$ can have
no heteroclinic intersections. This finishes the proof that $\{p_i\times A_i\}_{i=1}^k$ is a
skeleton.
\end{proof}

Let us proceed with the proof of Proposition~\ref{p.Kanseveralmeasures}.
Applying Proposition~\ref{p.newcretetion} to the elementary sets $\TT_i$ and the skeleton
$\{p_i\times A_i\}_{i=1}^k$ provided by Lemma~\ref{l.skeleton}, we find that $f$ has mostly
contracting center with $k$ physical measures $\mu_1, \dots, \mu_k$ such that $\supp\mu_i=\TT_i$
for every $i$.

\begin{lemma}\label{l.dense}
$W^s((p_i,A_i),f)$ is dense in $\TT^2\times S^2$ for every $i=1, \dots, k$.
\end{lemma}

\begin{proof}
By construction, the stable manifold of $A_i$ for the flow $X_i$ is dense in the sphere $W^c(p_i)$;
recall Figure~\ref{f.flow}.
It follows that the stable manifold $W^s((p_i,A_i),f)$ is dense in $W^s(W^c(p_i),f)$.
Moreover, the latter is dense in $\TT^2\times S^2$ because it coincides with $W^s(p_i,g)\times S^2$
and the stable manifold $W^s(p_i,g)$ is dense in $\TT^2$.
This proves the lemma.
\end{proof}

Then, by Theorem~\ref{t.mainA}, the basin of each physical measure $\mu_i$ is dense in $\TT^2\times S^2$.
This completes the proof of Proposition~\ref{p.Kanseveralmeasures} in what concerns the map $f$.
We are left to show that the conclusions extend to any $C^{1+\alpha}$ diffeomorphism $\tilde{f}$
in a neighborhood which leaves every $\TT_i$ fixed.

Begin by observing that $\tilde{f}\mid \TT_i$ is close to $f\mid \TT_i$ and, in particular, it is
Anosov. It follows that $\tilde{f}$ admits a unique Gibbs $u$-state supported on $\TT_i$
(the physical measure of that Anosov diffeomorphism) and that Gibbs $u$-state is close to $\mu_i$.
The latter ensures that the center Lyapunov exponents remain negative, and so $\TT_i$ remains an
elementary set for $\tilde{f}$. Each fixed point $(p_i,A_i)$ admits a continuation
$(p_i(\tilde{f}),A_i)$ for $\tilde{f}$. By Lemma~\ref{l.robustgraph}, these points form a pre-skeleton
for $\tilde{f}$. So, we are still in a position to use Proposition~\ref{p.newcretetion} to
conclude that $\tilde{f}$ has mostly contracting center and exactly $k$ physical measures,
$\tilde\mu_1, \dots, \tilde\mu_k$, with $\tilde\mu_i$ supported on $\TT^i$ for every $i$.
The proposition also states that  $\{p_i(\tilde f)\times A_i\}_{i=1}^k$ is actually a skeleton
for $\tilde{f}$.

We are left to prove that the basin of every $\tilde\mu_i$ is dense. By Theorem~\ref{t.mainA},
it suffices to show that the stable manifold of every $(p_i,A_i)$ is dense.
The center foliation of $f$ coincides with the trivial fibration $\{x\}\times S^2$,
which is normally hyperbolic and smooth. Thus, by the stability theorem of Hirsch, Pugh, Shub~\cite{HPS77},
the perturbation $\tilde{f}$ admits an invariant center foliation of $\TT^2\times S^2$ whose
leaves are $C^{1+\alpha}$ spheres uniformly close to the trivial fibers.
In particular, the center leaf through each point $(p_i(\tilde f),A_i)$ is close to
$\{p_i\} \times S^2$. That implies that the restriction of $\tilde f$ to that center leaf
is Morse-Smale and the stable manifold of $(p_i(\tilde f),A_i)$ is dense in it.
So, the stable manifold of $(p_i(\tilde{f}),A_i)$ is dense in the stable manifold of
$W^c(p_i(\tilde{f},A_i),\tilde{f})$.
The stability theorem also says that there exists a homeomorphism $h$ of $\TT^2\times S^2$ that
maps the center leaves of $f$ to the center leaves of $\tilde{f}$ and which is a leaf conjugacy:
$$
h(W^c(z,f))= W^c(h(z),\tilde{f})
\quad\text{for every } z\in \TT^2\times S^2.
$$
Then the stable manifold of $W^c(p_i(\tilde{f},A_i),\tilde{f})$ is just the image under $h$
of the stable manifold of $W^c(p_i,f)$. That guarantees that the stable manifold of
$W^c(p_i(\tilde{f},A_i),\tilde{f})$ is dense in $\TT^2 \times S^2$.
In this way we have recovered all the ingredients we used for $f$ and so at this point
our arguments extend to $\tilde{f}$, as claimed.

Finally, to construct examples with an odd number of physical measures, it suffices to show
that one can modify the diffeomorphism $f$ above, in such a way that the physical measures
of the resulting diffeomorphism $\hat f$ are precisely $\mu_1$, \dots, $\mu_{k-1}$.
Let $q\in T^2$ be a point of transverse intersection between $W^u(p_k,g)$ and $W^s(p_1,g)$.
Let $Y^t$, $t\in\RR$ be a smooth flow on $T^2\times S^2$ such that 
\begin{enumerate}
\item $Y^t$ is supported on a small neighborhood of $(q,A_k)$;
\item $Y^t$ preserves the center foliation;
\item for any $t>0$, the map $Y^t$ sends $(q,A_k)$ to some $(q,C)$ with $C\notin\cC$.
\end{enumerate}
Pick $\hat f = Y^t \circ f$ for any $t>0$. Condition (3) implies that for $\hat f$ the
unstable manifold of $(p_k,A_k)$ intersects the stable manifold of $(p_1,A_1)$. So 
\begin{equation}\label{eq.skel}
\{(p_1, A_1), \dots, (p_{k-1},A_{k-1})\}
\end{equation}
is a pre-skeleton for $\hat f$. Conditions (1) and (2) ensure that the unstable manifold of
each $(p_i,A_i)$ remains unperturbed and thus is still contained in $T^2 \times\{A_i\}$,
for $i=1, \dots, k-1$. This ensures that the set in \eqref{eq.skel} is actually a skeleton,
and so $\hat f$ has exactly $k-1$ physical measures. All the other stated properties are
obtained just as in the previous case.
\end{proof}

\subsection{Collapse of measures and explosion of basins\label{ss.examples}}

In this subsection, we prove that the examples we have just constructed are statistically \emph{unstable}:
the simplex generated by all the physical measures does not vary continuously with the dynamics,
as physical measures may collapse, with their basins of attraction exploding, after small perturbations of
the diffeomorphism. In fact, we obtain two different instability results:

\begin{itemize}
\item For any proper subset of physical measures, one can find a small perturbation of the
original diffeomorphism for which those physical measures vanish: their basins are engulfed
by the ones of the remaining physical measures.

\item For any Gibbs-$u$ state $\mu$ of the original diffeomorphism (not necessarily ergodic),
one can find diffeomorphisms $f_n$ converging to $f$, such that every each $f_n$ has a unique
physical measure $\mu_n$ and the sequence $(\mu_n)_n$ converges to $\mu$ in the weak-* topology.
\end{itemize}

In all that follows $f:M\to M$ is a partially hyperbolic diffeomorphism with $k=3$ physical measures,
as constructed in the previous section (the constructions extend to arbitrary $k$ in a straightforward way).
Let us first describe our perturbation technique.
It is designed to create new heteroclinic intersections, thus reducing the number of saddle points
in the skeleton.

For distinct $i, j \in\{1, 2, 3\}$, let $q_{i,j}\in T^2$ be a point of transverse intersection of
$W^u(p_i,g)$ and $W^s(p_j,g)$. Consider a smooth flow $Y_{i,j}^t$ on $T^2\times S^2$ such that:
\begin{enumerate}
\item $Y_{i,j}^t$ is supported on a small neighborhood of $(q_{i,j}, A_i)$;
\item $Y_{i,j}^t$ preserves the center foliation of $f$;
\item for any $t>0$, the map $Y_{i,j}^t$ sends $(q_{i,j}, A_i)$ to some $(q_{i,j}, C_i)$
with $C_i\notin \cC$.
\end{enumerate}
We will always consider perturbations $f_{t_1, t_2, t_3}$ of the original $f$ of the form
$$
f_{t_1,t_2,t_3} = Y_{1,2}^{t_1} \circ Y_{2,3}^{t_2} \circ Y_ {3,1}^{t_3} \circ f
\quad\text{$t_1$, $t_2$, $t_3$ close to zero.}
$$
Observe that $f_{t_1,t_2,t_3}| p_i\times S^2=f|p_i\times S^2$, since $p_i\times S^2$, $i=1,2,3$
are away from the regions of perturbation. By Lemma~\ref{l.robustgraph},
$\{p_i\times A_i\}_1^3$ is a pre-skeleton of $f_{t_1,t_2,t_3}$.
Denote $p_4=p_1$ and $A_4=A_1$ and $q_{3,4}=q_{3,1}$.

\begin{lemma}\label{l.perturbation}
The strong unstable leaf $\cF^u((p_i, A_i), f_{t_1,t_2,t_3})$ has some transverse intersection
with $W^s((p_{i+1},A_{i+1}),f_{t_1,t_2,t_3})$, for every $t_i>0$.
\end{lemma}

\begin{proof}
Let $j=i+1$. By construction, the strong unstable leaf of $(p_i,A_i)$ for $t_{t_1, t_2, t_3}$
contains the point $Y_{i,j}^{t_i}(q_{i,j},A_i)$, which is the strong stable leaf of some point in
$\{p_j\} \times (S^2\setminus\cC)$. The latter is in the stable manifold of $(p_j,A_j)$.
Clearly, the two manifolds intersect transversely at this point.
\end{proof}

We are ready to state and prove our first instability result:

\begin{proposition}\label{p.weakcollapse}
Given any proper subset $\Gamma$  of the set $\{\mu_1, \mu_2, \mu_3\}$ of physical measures of $f$,
one can find $\tilde{f}$ arbitrarily close to $f$ such that the set of physical measures of $\tilde{f}$
is $\{\mu_1,\mu_2,\mu_3\}\setminus \Gamma$.
\end{proposition}

\begin{proof}
First, suppose that $\#\Gamma=1$, say, $\Gamma=\{\mu_1\}$. Consider $\tilde{f}=f_{t_1,0,0}$ with $t_1>0$.
The measures $\mu_2$ and $\mu_3$ are still ergodic Gibbs-$u$ states and physical measures for $\tilde{f}$,
since $\tilde{f}$ coincides with $f$ on the neighborhood of their supports, $T_2$ and $T_3$.
Moreover, the unstable manifolds of $(p_2, A_2)$ and $(p_3, A_3)$ are still contained in $T_2$ and
$T_3$, respectively, and so these points have no heteroclinic intersections.
On the other hand, by Lemma~\ref{l.perturbation}, $(p_1, A_1) \prec (p_2, A_2)$.
Thus, $\{(p_2, A_2), (p_3, A_3)\}$ is a skeleton of $\tilde{f}$, by Lemma~\ref{l.existenceskeleton}.
So, by Theorem~\ref{t.mainA}, the diffeomorphism $\tilde{f}$ has exactly two physical measures,
$\mu_1$ and $\mu_2$.

Now suppose that $\#\Gamma=2$, say, $\Gamma=\{\mu_1,\mu_2\}$. Consider $\tilde{f}=f_{t_1,t_2,0}$
with $t_1>0$ and $t_2>0$. Then, just as before, $(p_1, A_1) \prec (p_2, A_2)$ and
($p_2, A_2) \prec (p_3, A_3)$ for $\tilde{f}$. Then, by Lemma~\ref{l.existenceskeleton},
$\{(p_3,A_3)\}$ is a skeleton of $\tilde{f}$ and, by Theorem~\ref{t.mainA}, the map
$\tilde{f}$ has a unique physical measure, $\mu_3$.
\end{proof}

The same arguments show that if $t_1, t_2, t_3$ are all positive then $f_{t_1, t_2, t_3}$
has a unique physical measure (the points $(p_i, A_i)$ are all heteroclinically related),
which need not be close to any of the physical measures of the original map $f$.

\begin{proposition}\label{p.strongcollapse}
For each Gibbs $u$-state $\nu$ of $f$ there exists a sequence $(f_n)_n$ such that every $f_n$
has a unique physical measure $\mu_n$ and the sequence $(\mu_n)_n$ converges to $\nu$ as $n\to\infty$.
\end{proposition}

\begin{proof}
Notice that $\nu$ is an element of the simplex
$$
\Delta=\{(s_1\mu_1+s_2\mu_2+s_3\mu_3): s_1\geq 0, s_2\geq 0, s_3\geq 0, s_1+s_2+s_3=1\},
$$
since every Gibbs $u$-state is a linear combination of the ergodic Gibbs $u$-states and,
for $f$, these are precisely the physical measures. Clearly, it is no restriction to suppose
that $\mu$ belongs to the interior of $\Delta$. Let $P:\cM\to V$ be any continuous affine map
from the Banach space $\cM$ of finite signed measures on $\TT^2\times S^2$ to the affine plane
$V\subset\cM$ generated by $\Delta$ such that $P \mid V = \id$. Existence of such a map follows
from the Hahn-Banach theorem.

For $n\ge 1$ and $0<\delta<1/n$, consider the hexagon
$$
\begin{aligned}
H_n(\delta) = \{(t_1,t_2,t_3): t_1\geq 0, t_2\geq 0, t_3\geq 0, t_1+t_2+t_3=\frac 1n,\\
t_1+t_2\ge\delta, t_1+t_3\ge\delta, t_2+t_3\ge\delta \}.
\end{aligned}
$$
Every triple $(t_1,t_2,t_3)\in H_n(\delta)$ has at least two positive coordinates. Hence, by the
same arguments as in the proof of Proposition~\ref{p.weakcollapse}, the corresponding map
$f_{t_1, t_2, t_3}$ has exactly one Gibbs $u$-state $\mu_{t_1,t_2,t_3}$, which is also the unique
physical measure. This defines a map $\Phi(t_1,t_2,t_3)=\mu_{t_1,t_2,t_3}$ with values in the
space of probability measures on $\TT^2\times S^2$. By upper semi-continuity of the space of Gibbs
$u$-states, $\Phi$ is continuous on $H_n(\delta)$ and the image $\Phi(H_n(\delta))$ is contained
in a neighborhood of the simplex $\Delta$.

Let $\alpha$ be the distance from $\mu$ to the boundary of $\Delta$. We claim that for each $n$
there exists $0 < \delta_n < 1/n$ such that the image of $\tilde H_n=H_n(\delta_n)$ under $\Phi$
is a topological simplex $(\alpha/4)$-close to $\Delta$ in the space $\cM$, in the following sense:
\begin{itemize}
\item[(i)] the two simplices have the same vertices and
\item[(ii)] every edge of $\Phi(\tilde H_n)$ is contained in the $(\alpha/4)$-neighborhood of the
corresponding edge of $\Delta$.
\end{itemize}\
It follows that for every large $n$ the image $P(\Phi(\tilde H_n))$ is a topological simplex
$(\alpha/2)$-close to $\Delta$ in the plane $V$. By a topological degree argument, it follows that
$P(\Phi(\tilde H_n))$ contains $\mu$: otherwise, it would be retractable to the boundary of $\Delta$,
which is nonsense. This means that there exists $(t_1(n),t_2(n),t_3(n))\in\tilde H_n$ such that
$P(\mu_{t_1(n),t_2(n),t_3(n)})=\mu$. Let
$$
f_n=f_{t_1(n), t_2(n), t_3(n)} \quand \mu_n=\mu_{t_1(n), t_2(n), t_3(n)}.
$$
The definition of $\tilde H_n$ implies that $t_i(n)\to 0$ when $n\to \infty$ for every $i$. Thus,
$(f_n)_n$ converges to $f$. By upper semi-continuity of the space of Gibbs $u$-states, every
accumulation point of the sequence $(\mu_n)_n$ is contained in $\Delta$. Also, by construction,
$P(\mu_n)=\mu$ for every $n$. Since $P$ is  continuous and its restriction to $\Delta$ is
injective, this implies that $(\mu_n)_n$ converges to $\mu$.

We are left to prove the claim above. Let $I_1, I_2, I_3$ and
$J_1, J_2, J_3$ be the boundary segments of $H_n(\delta)$, with $I_j$ contained in $\{t_j=0\}$
and $J_j$ contained in $\{t_j+t_{j+1}=\delta\}$ (denote $t_4=t_1$ and $\mu_4=\mu_1$).
If $(t_1,t_2,t_3)\in I_j$ then $t_j$ is the unique vanishing parameter and so
$\mu_{t_1,t_2,t_3}=\mu_j$. This means that $\Phi(I_j)=\{\mu_j\}$ for  $j=1, 2, 3$, which
gives part (i) of the claim. It also follows that $\Phi(J_j)$ is a continuous curve from
$\mu_j$ to $\mu_{j+1}$. Using upper semi-continuity once more, this curve must be contained
in the $(\alpha/4)$-neighborhood of the space of Gibbs $u$-states of $f_{t_1,t_2,t_3}$
with $t_j=t_{j+1}=0$, provided $\delta$ is small enough. To conclude, it suffices to observe
that the latter is precisely the edge $[\mu_j,\mu_{j+1}]$ of $\Delta$.
\end{proof}

\section{Continuity of basins}

In this section, we prove Theorem~\ref{t.mainC}. Indeed, we prove the following somewhat more
explicit fact:

\begin{proposition}\label{p.continuousbasin}
Let $k\ge 1$ and $\cO$ be a subset of the space of $C^{1+\alpha}$ diffeomorphisms of $M$ with
mostly contracting center such that every $f\in\cO$ has exactly $k$ physical measures,
$\mu_1(f)$, \dots, $\mu_k(f)$. Let $\{f_n\}_{n=1}^\infty$ be any sequence in $\cO$ converging
to some $f\in\cO$.  Then, up to suitable numbering,
$$
d\big(B(f_n,\mu_i(f_n)), B(f,\mu_i(f))\big) \to 0 \quad\text{for every $1\leq i\leq k$.}
$$
\end{proposition}

The conclusion holds, in particular, within the family of examples constructed in
Subsection~\ref{Kan's example} or, more precisely, in the last part of
Proposition~\ref{p.Kanseveralmeasures}.

\begin{proof}
Let $\{p_1(f),\cdots, p_k(f)\}$ be a skeleton of $f\in\cO$ with $p_i(f)\in \supp \mu_i(f)$ for each $i$.
As we have seen before, the continuations $p_i(g)$, $1 \leq i \leq k$ of the saddle points $p_i(f)$
constitute a skeleton for every $g$ in a small neighborhood relative $\cU\subset\cO$ (because the
number of physical measures remains the same).

\begin{lemma}\label{l.ref1}
Let $l$ be the product of the periods of $p_1(f)$, \dots, $p_k(f)$. If $\nu$ is a Gibbs $u$-state
of any iterate $f^n$, $n\ge 1$ then
\begin{equation}\label{eq.Gibbssum}
\frac 1l \sum_{j=0}^{l-1} f_*^j\nu
\text{ is a Gibbs $u$-state of $f$.}
\end{equation}
\end{lemma}

\begin{proof}
We begin by claiming that $f^l_*\nu=\nu$. For proving this claim, it suffices to consider the case
when $\nu$ is ergodic for $f^n$. Notice that $f^n$ has mostly contracting center and
$$
\{f^j(p_i(f)): 1 \le i \ k \text{ and } 0 \le j < \per(p_i(f))\}
$$
is a pre-skeleton for $f^n$. Thus, by Theorem~\ref{t.mainA}, the support of $\nu$ contains some
$f^j(p_i)$. The measure $f^l_*\nu$ is still $f^n$-invariant and $f^n$-ergodic. Then, since $f$
preserves absolute continuity along unstable manifolds, $f_*^l\nu$ is still a Gibbs $u$-state
for $f^n$. Since its support also contains $f^j(p_i)$, it follows from Theorem~\ref{t.mainA} that
$f_*^l\nu$ and $\nu$ coincide, as claimed. Then the measure in \eqref{eq.Gibbssum} is $f$-invariant
and, using once more the fact that $f$ preserves absolute continuity along unstable manifolds,
it is a Gibbs $u$-state for $f$, as we wanted to prove.
\end{proof}

\begin{lemma}\label{l.ref2}
There exists $\lambda_0>0$ and for every large $N\ge 1$ there exists a relative neighborhood
$\cU_N\subset\cO$ of $f$ such that
$$
\int \log \|Dg^N \mid E^{cs}\| \, d\nu \le - 2 N \lambda_0
$$
for every $g\in\cU_N$ and any Gibbs $u$-state $\nu$ of $g^N$.
\end{lemma}

\begin{proof}
Since $f$ has mostly contracting center, the largest center exponent
$$
\lim_n \int \frac{1}{n} \log \|Df^n \mid E^{cs}\| \, d\mu_i,
$$
is negative for every $i=1, \dots, k$. Since every Gibbs $u$-state is a convex
combination of $\mu_1(f)$, \dots, $\mu_k(f)$, it follows that there exist $N\ge 1$
and $\lambda_0>0$ such that
\begin{equation}\label{eq.Gibs1}
\int \log \|Df^N \mid E^{cs}\| \, d\mu \le - 8 N \lambda_0
\end{equation}
for every Gibbs $u$-state $\mu$ of $f$.
Now let $\nu$ be any Gibbs $u$-state for $f^N$. It follows from \eqref{eq.Gibs1} and
Lemma~\ref{l.ref1} that
\begin{equation}\label{eq.Gibs2}
\frac 1l \sum_{j=0}^{l-1} \int \log \|Df^N \mid E^{cs}\|\circ f^j \, d\nu \le - 4 N \lambda_0.
\end{equation}
For any $N\ge l$ and $x\in M$, we have
$$
Df^N \mid E^{cs}(f^j(x)) = Df^j \mid E^{cs}(f^N(x)) \circ Df^N \mid E^{cs}(x) \circ Df^{-j} \mid E^{cs}(f^j(x)).
$$
Hence, denoting $C=\max \log\|Df\| + \max \log\|Df^{-1}\|$,
$$
\log \|Df^N \mid E^{cs}\| \circ f^j
\ge \log \|Df^N \mid E^{cs}\| - C j
\ge \log \|Df^N \mid E^{cs}\| - C l.
$$
Combining this inequality with \eqref{eq.Gibs2}, we obtain
$$
\int \log \|Df^N \mid E^{cs}\| \, d\nu \le - 4 N \lambda_0 + C l \le - 3 N \lambda_0
$$
as long as we take $N \ge C l/\lambda_0$.
By upper semi-continuity of the set of Gibbs $u$-states, it follows that
$$
\int \log \|Dg^N \mid E^{cs}\| \, d\nu \le - 2 N \lambda_0
$$
for any Gibbs $u$-state $\nu$ of $g^N$ and any $g$ in a neighborhood $\cU_N$ of $f$.
\end{proof}

Fix $N$ to be a multiple of $l$ large enough that Lemma~\ref{l.ref2} is satisfied.
For each $1\leq i \leq k$, choose a small neighborhood  $V_i$ of $p_i(f)$. Fix $\rho>0$ small,
such that the $\rho$-neighborhood $W^u_\rho(p_i(g),g)$ of $p_i(g)$ inside its unstable manifold
$W^u(p_i(g),g)$ is contained in $V_i$ for every $g\in\cU$ and every $i$.
For each $g\in\cU$ and $i=1, \dots, k$, define $\Lambda_i(m,g)$ to be the subset of points
$x\in W^u_{2\rho}(p_i(g),g)$ such that
$$
\frac 1n \sum_{j=0}^{n-1} \log \|Dg^{N} \mid E^{cs}(g^{jN}(x))\| \leq - N\lambda_0 \text{ for all } n \geq m.
$$

%
%

\begin{lemma}\label{l.ref3}
There exists a relative neighborhood $\cU\subset\cU_N$ of $f$ and there exist constants $K>0$ and $\theta<1$ such that
$$
\Leb(W^u_{2\rho}(p_i(g),g) \setminus \Lambda_i(m,g)) \le K\theta^m \quad\text{for every $g\in\cU$ and $m\ge 1$.}
$$
\end{lemma}

\begin{proof}
Let $g\in\cU_N$ and $l$ be the Lebesgue measure on the $u$-disk $W^u_\rho(p_i(g),g)$.
Define
$$
A(x) = \|Dg^{N} \mid E^{cs}(g^{N}(x))\|
$$
for $x\in M$. Then let $I$ be the set (compact interval) of values of $\int A \, d\nu$ over all Gibbs $u$-states of $g^N$.
By Lemma~\ref{l.ref2}, $I$ is contained in $(-\infty,-2N\lambda_0)$. Then, for each fixed $g$, the claim is contained in
the conclusion of~\cite[Proposition~1]{Dol_un} (or \cite[Theorem~1]{Dol04a}, in the special case when there exists a
unique Gibbs $u$-state) for $\vep=\lambda_0$). Moreover, the constants may be taken uniform over all $g$ in some
neighborhood $\cU\subset\cU_N$ of $f$: see \cite[Section~7]{Dol04a} for the case when there is a unique Gibbs $u$-state,
and \cite[Exercise~7]{Dol_un} for the general case.
\end{proof}

For each large $m\ge 1$ there exists $\delta_m>0$ such that for any $g\in\cU$ the Pesin stable manifold of every
$y \in \Lambda_i(m,g)$ has uniform size $\delta_m$ (meaning that it contains a $\dim E^{cs}$-disk of radius
$\delta_m$ around $y$). Indeed, the uniform bound on the size of the stable manifold follows from the same arguments as
\cite[Lemma~3.7]{ABV00}, applied to the inverse of $g$.

Moreover, these local Pesin stable manifolds define a lamination $W^s(\Lambda_{i}(m,g))$ which is absolutely
continuous (see \cite{Pes76,PSh89}): the corresponding holonomy maps $h:D_1 \to D_2$ between disks $D_1$ and
$D_2$ transverse to the lamination are absolutely continuous, with Jacobian given by
\begin{equation}
J h(y)
= \lim_k \frac{\det(Dg^k \mid T_y D_1)}{\det(Dg^k \mid T_{h(y)} D_2)}
= \prod_{j=0}^\infty \frac{\det(Dg \mid T_{g^j(y)} g^j(D_1))}{\det(Dg \mid T_{g^j h(y)} g^j(D_2))}
\end{equation}
for $y\in D_1 \cap W^s(\Lambda_{i}(m,g))$. In particular, for any $m\ge 1$ there exists $\gamma_m>0$ such that the
Jacobian is bounded above by $2$, for any $g\in \cU$ and any disks $D_1$ and $D_2$ in the $\gamma_m$-neighborhood of
$W^u_{2\rho}(p_i(g),g)$ in the $C^1$ topology.

Let $K>1$ be an upper bound for the distortion of backward iterates of any $f$ along unstable disks:
\begin{equation}\label{eq.distortion_unstable}
\frac{\det(Df^{-n} \mid T_{x_1} D)}{\det(Df^{-n} \mid T_{x_2} D)} \le K
\end{equation}
for any $x_1, x_2 \in D$ and any $u$-disk of $f$ with radius $1$.
Fix $\kappa>0$ such that $\Leb(W^u_{\rho/2}(z))\ge\kappa$ for every $z\in M$.

\begin{lemma}\label{l.fixm}
Given $\vep>0$ there exists $m\ge 1$ such that
$$
\Leb\big(D \setminus B(g,\mu_i(g)))\big) < \frac{\vep\kappa}{4K}
$$
for any $g\in \cU$ and any disk $D$ in the $\gamma_m$-neighborhood of $W^u_{2\rho}(p_i(g),g)$ in the $C^1$-topology.
\end{lemma}

\begin{proof}
By~\cite[Theorem~11.16]{Beyond}, Lebesgue almost every point in $W^u(p_i(g),g)$ is in the basin of
some Gibbs $u$-state. Since $p_i(g)$ is in the support of $\mu_i(g)$, by the definition of skeleton,
and the supports are disjoint, we get that almost every point in $W_{2\rho}^u(p_i(g),g)$ is in
$B(g,\mu_i(g))$. Then the same is true for (every point in the Pesin stable manifold through)
almost every point in $\Lambda_i(m,g)$. By Lemma~\ref{l.ref3}, we may fix $m\ge 1$ such that the
Lebesgue measure of the complement of $\Lambda_i(m,g)$ in $W^u_{2\rho}(p_i(g),g)$ is less than
$(\vep\kappa)/(8K)$. In view of the previous observations, and the fact that the Jacobian is bounded by $2$,
it follows that the Lebesgue measure of the complement of $B(\mu_i(g))$ in $D$ is less than
$(\vep\kappa)/(4K)$ for any $D$ in the $\gamma_m$-neighborhood of $W^u_{2\rho}(p_i(g),g)$, as claimed.
\end{proof}

Now we apply to the diffeomorphism $f$ the local Markov construction in \cite[Section~4.2]{Almost}:
for any small $\delta>0$ we may find a family $\{U(z): z \in W^s_\delta(p_i,f)\}$ of
embedded $u$-disks such that
$$
W^u_\rho(z) \subset U(z) \subset W^u_{2\rho}(z)
\quad\text{for every $z\in W^s_\delta(p_i)$,}
$$
and, for any $j\ge 1$ and $z, \zeta \in W^s_\delta(p_i)$, either
$$
f^{-j}(U(z)) \cap U(\zeta) = \emptyset \quad\text{or} \quad f^{-j}(U(z))\subset U(\zeta).
$$

Given $\vep>0$, fix $m\ge 1$ as in Lemma~\ref{l.fixm} from now on. Then, take $\delta\in(0,\delta_m)$
such that $W^u_{2\rho}(z)$ is in the $(\gamma_m/2)$-neighborhood of $W^u_{2\rho}(p_i)$ for every
$z \in W^s_\delta(p_i)$. Denote by $\cF_i$ the union of the disks $U(z)$, $z\in W^s_\delta(p_i)$.

\begin{lemma}\label{l.fillin_forf}
Up to zero Lebesgue measure,
$$
\bigcup_{j=0}^\infty f^{-j}\big(\cF_i \cap W^s(\Lambda_{i}(m,f))\big) = B(\mu_i).
$$
\end{lemma}

\begin{proof}
By~\cite[Theorem~11.16]{Beyond}, Lebesgue almost every point in the unstable manifold of $p_i$ is in the basin of
some Gibbs $u$-state. Recalling that $p_i$ is in the support of $\mu_i$, by the definition of skeleton,
and the supports are disjoint, we get that almost every point in $W_{2\rho}^u(p_i)$ is in the basin $B(\mu_i)$.
Since the basin is saturated by stable sets, and the lamination
$\cW^s(\Lambda_{i}(m,f)$ is absolutely continuous, it follows that
$$
W^s(\Lambda_{i}(m,f)) \subset B(\mu_i)
\quad\text{up to zero Lebesgue measure.}
$$
Since the basin is an invariant set, this implies the inclusion $\subset$ in the statement.

The converse is a corollary of \cite[Proposition~6.9]{ViY13}. Indeed, this proposition implies that
$\cV = \cup_{i=1}^k\cup_{j=0}^\infty f^{-j}(\cF_i \cap W^s(\Lambda_{i}(m,f)))$
contains a full Lebesgue measure subset of every strong-unstable disk. By the absolute continuity of the
strong unstable foliation, this implies that $\cV$ contains a full volume subset of the ambient manifold.
Since we already know that each $\cV_j=\cup_{j=0}^\infty f^{-j}(\cF_i \cap W^s(\Lambda_{i}(m,f)))$
is contained in the corresponding basin $B(\mu_j)$, and the basins are pairwise disjoint,
it follows that $\cV_j = B(\mu_j)$ up to measure zero. The proof is complete.
\end{proof}

By Lemma~\ref{l.fillin_forf}, we may fix $N\ge 1$ such that
\begin{equation}\label{eq.measure1}
\Leb\big(B(\mu_i) \setminus \bigcup_{j=0}^N f^{-j}(\cF_i)\big) < \frac{\vep}{2}.
\end{equation}
In view of our choice of $\delta>0$, we may find a neighborhood $\tilde\cU\subset\cU$ of $f$ such that
$g^jf^{-j}(U(z))$ is contained in some disk $D$ in the $\gamma_m$-neighborhood of $W^u_{2\rho}(p_i(g),g)$
for every $z\in W^s_\delta(p_i)$ and $0 \le j \le N$. Reducing $\tilde\cU$ if necessary, and recalling
\eqref{eq.distortion_unstable}, we may suppose that
$$
\frac{\det(Dg^{-j} \mid T_{x_1} D)}{\det(Dg^{-j} \mid T_{x_2} D)} \le 2K
$$
for any $0 \le j \le N$, any $x_1, x_2 \in D$ and any disk $D$ in the $\gamma_m$-neighborhood of
$W^u_{2\rho}(p_i(g),g)$. It is clear that
$g^jf^{-j}(U(z))$ converges to $U(z) \supset W_\rho^u(z)$ when $^g \to f$. Thus, recalling our choice of
$\kappa$ and further reducing $\tilde\cU$ if necessary, we may suppose that
$$
\Leb\big(g^jf^{-j}(U(z))\big) \ge \kappa
$$
for any $0 \le j \le N$ and any $g\in\tilde\cU$.

\begin{lemma}\label{l.bounded_distortion_basin}
For every $z\in W_\delta^s(p_i)$ and $j \ge 0$ and $g\in\tilde\cU$,
$$
\Leb\big(f^{-j}(U(z)) \setminus B(g,\mu_i(g)))\big) < \frac{\vep}{2} \Leb\big(f^{-j}(U(z))\big).
$$
\end{lemma}

\begin{proof}
Let $D$ be a disk in the $\gamma_m$-neighborhood of $W^u_{2\rho}(p_i(g),g)$ and containing $g^jf^{-j}(U(z))$.
By Lemma~\ref{l.fixm},
$$
\Leb\big(g^jf^{-j}(U(z)) \setminus B(g,\mu_i(g))\big)
\le \Leb\big(D \setminus B(g,\mu_i(g))\big) < \frac{\vep\kappa}{4K}
$$
and so,
$$
\frac{\Leb\big(g^jf^{-j}(U(z)) \setminus B(g,\mu_i(g))\big)}{\Leb\big(g^jf^{-j}(U(z))}
< \frac{\vep}{4K}.
$$
Then, since the basin $B(g,\mu_i(g))$ is a $g$-invariant set,
$$
\frac{\Leb(f^{-j}(U(z)) \setminus B(\mu_i))}{\Leb(f^{-j}(U(z)))}
\le 2 K \frac{\Leb\big(g^jf^{-j}(U(z)) \setminus B(g,\mu_i(g))\big)}{\Leb\big(g^jf^{-j}(U(z))}
< \frac{\vep}{2},
$$
as claimed.
\end{proof}

\begin{corollary}\label{c.foliatedbasin}
For every $g\in\tilde\cU$,
\begin{equation}\label{eq.measure2}
\Leb\big(\bigcup_{j=0}^N f^{-j}(\cF_i) \setminus B(g,\mu_i(g))\big) < \frac{\vep}{2}.
\end{equation}
\end{corollary}

\begin{proof}
Define the \emph{return time} $r(z)$ of each $z \in W^s_\delta(p_i)$ to be the smallest $n \in\NN\cup\{\infty\}$
such that $f^{n}(U(z))$ intersects (and thus is contained in) some $U(\zeta)$, $\zeta\in W_\delta^s(p_i)$. Observe that
$$
\bigcup_{j=0}^N f^{-j}(\cF_i)
$$
is the pairwise disjoint union of the pre-images $f^{-j}(U(z))$ with $z\in W^s_\delta(p_i)$ and
$0 \le j \le \min{r(z)-1,N}$. For each one of these pre-images, Lemma~\ref{l.bounded_distortion_basin} gives that
$$
\Leb\big(f^{-j}(U(z)) \setminus B(g,\mu_i(g)))\big) < \frac{\vep}{2} \Leb\big(f^{-j}(U(z))\big).
$$
So, by the Cavalieri principle,
$$
\Leb\big(\bigcup_{j=0}^N f^{-j}(\cF_i) \setminus B(g,\mu_i(g))\big)
< \vep \Leb\big(\bigcup_{j=0}^N f^{-j}(\cF_i)\big) \le \frac{\vep}{2}.
$$
This proves the claim.
\end{proof}

Combining \eqref{eq.measure1} and \eqref{eq.measure2}, we find that
\begin{equation}\label{eq.measure3}
\Leb\big(B(\mu_i) \setminus B(g,\mu_i(g))\big) < \vep
\quad\text{for every $g \in \tilde\cU$.}
\end{equation}
Since, for both $f$ and $g$, the basins are pairwise disjoint and their union has total measure,
$$
B(g,\mu_i(g))\setminus B(\mu_i) \subset \bigcup_{j\neq i} B(\mu_j) \setminus B(g,\mu_i(g))
$$
up to measure zero, for every $i$. Thus, it also follows from \eqref{eq.measure3} that
\begin{equation}\label{eq.measure4}
\Leb\big(B(g,\mu_i(g))\setminus B(\mu_i)\big) < (k-1) \vep
\quad\text{for every $g \in \tilde\cU$.}
\end{equation}
The relations \eqref{eq.measure3} and \eqref{eq.measure4} mean that $d(B(\mu_i),B(g,\mu_i(g)) < k \vep$
for every $g\in\tilde\cU$, and so the argument is complete.
\end{proof}


\end{document}